\documentclass[12pt]{article}
\usepackage{hyperref}
\usepackage{color}
\hypersetup{colorlinks=true}
\usepackage{a4wide}
\usepackage{amsmath}
\usepackage{amssymb}
\usepackage{mathrsfs}
\usepackage[english,francais]{babel}
\usepackage[fixlanguage]{babelbib}
\usepackage{theorem}
\usepackage{amsfonts}
\usepackage[all]{xy}
\btxprintISSN{false}
\btxprintISBN{false}
\newlength{\secskip}
\newlength{\thsecskip}
\newlength{\ssecskip}
\newlength{\thssecskip}
\newlength{\sssecskip}
\newlength{\thsssecskip}
\makeatletter
\renewcommand\section{\@startsection{section}{1}{\z@}%
                                   {\secskip}
                                   {2.3ex \@plus.2ex}
                                   {\normalfont\Large\bfseries}}
\renewcommand\subsection{
			  \@startsection{subsection}{2}{\z@}%
                                     {\ssecskip}
                                     {1.5ex \@plus .2ex}
									 {\normalfont\normalsize\bfseries}}
\renewcommand\subsubsection{
              \@startsection{subsubsection}{3}{\z@}%
                                     {\sssecskip}
                                     {\z@}
									 {\normalfont\normalsize\bfseries}}
\makeatother

\makeatletter
\renewcommand\l@subsection{\@dottedtocline{2}{1.5em}{2.7em}}
\makeatother

\makeatletter
\renewcommand\l@subsubsection{\@dottedtocline{3}{3.8em}{3.5em}}
\makeatother
\newcommand\rref[1]{{\rm\ref{#1}}}
\theoremstyle{change}
\theorembodyfont{\slshape}
\newtheorem{thm}[subsection]{Th\'eor\`eme.}

\newtheorem{subthm}[subsubsection]{Th\'eor\`eme.}
\newtheorem{subprop}[subsubsection]{Proposition.}
\newtheorem{subcor}[subsubsection]{Corollaire.}

\theorembodyfont{\rmfamily}

\newtheorem{xrems}[subsection]{Remarques.}

\newtheorem{subrem}[subsubsection]{Remarque.}

\newcounter{nr}
\renewcommand{\thenr}{{\bf\thesubsection.\arabic{nr}}}
\newenvironment{rems}%
{\begin{xrems}{\ }\par
\begin{list}{\thenr}{\usecounter{nr}\setlength{\leftmargin}{0cm}\setlength{\labelwidth}{2em}\setlength{\labelsep}{.6em}\setlength{\itemindent}{2.6em}}
}%
{\end{list}%
\end{xrems}}
\newcounter{nc}
\renewcommand{\thenc}{{\rm(\roman{nc})}}
\newenvironment{romlist}%
{\begin{list}{\thenc}{
\usecounter{nc}
\parsep=0pt
\setlength  \labelwidth{\leftmargin}
\addtolength\labelwidth{-\labelsep}
}
}{\end{list}} 
\newcommand{\Spec}{{\rm Spec\,}}
\newcommand{\pr}{{\rm pr}}

\newcommand{\Aa}{\mathbb{A}} 
 
\newcommand{\RR}{\mathbb{R}} 
\newcommand{\wh}[1]{\widehat{#1}}
\newcommand{\carrenoir}{\vrule height2.5mm width2mm depth 0mm}
\def\qed{ \hbox to5mm{}\hfill\llap{\carrenoir}}
 
\newcommand{\dem}{{\sl \noindent D\'emonstration\hskip0.1em: }} 
\newcommand{\inj}{\hookrightarrow}

\newcommand{\cO}{\mathscr{O}}
\newcommand{\cF}{\mathscr{F}}
\newcommand{\Top}{\mathrm{top}}
\date{}
\title{Un th\'eor\`eme de l'application ouverte sur les corps valu\'es alg\'ebriquement clos}
\author{Laurent Moret-Bailly
\thanks{L'auteur est membre du projet \og Points entiers et points rationnels\fg\ de l'Agence nationale pour la recherche.}
\medskip
\\
{\small IRMAR (Institut de Recherche Math\'{e}matique de Rennes,} \\
{\small UMR 6625 du CNRS)}\\
{\small Universit\'{e} de Rennes 1}, 
{\small Campus de Beaulieu, F-35042 Rennes Cedex}\\
{\small \href{mailto:laurent.moret-bailly@univ-rennes1.fr}{laurent.moret-bailly@univ-rennes1.fr}}\\
{\small{\href{http://perso.univ-rennes1.fr/laurent.moret-bailly/}{http://perso.univ-rennes1.fr/laurent.moret-bailly/}}}\bigskip\\
\`A para\^{\i}tre dans \href{http://www.mscand.dk/}{\textsl{Mathematica Scandinavica}}\\
\small{(accept\'e le 31 mai 2011)}}
\begin{document}
\selectlanguage{french}
\maketitle
\begin{abstract}
Soit $K$ un corps alg\'ebriquement clos muni d'une valeur absolue non triviale, et soit $f:X\to Y$ un morphisme de $K$-sch\'emas de type fini. On montre que $f$ est  uni\-ver\-sel\-le\-ment ouvert si et seulement si l'application induite sur les points $K$-rationnels est ouverte pour les topologies d\'eduites de la valeur absolue.
\end{abstract}\selectlanguage{english}
\begin{abstract}
Let $K$ be an algebraically closed field with a nontrivial absolute value, and let  $f:X\to Y$ be a morphism of $K$-schemes of finite type. We show that $f$ is  universally open if and only if  the induced map on $K$-rational points is open for the topologies deduced from the absolute value.
\end{abstract}
\selectlanguage{french}
\noindent\textsl{Classification AMS:\/} 14B25, 12J10.
\section{Introduction}
Soit $(K,\vert\,.\,\vert)$ un corps valu\'e (c'est-\`a-dire que $\vert\,.\,\vert: K\to\RR$ est une valeur absolue sur $K$, qui sera suppos\'ee non triviale).  Pour tout $K$-sch\'ema de type fini $X$, on notera $X_{\Top}$ l'ensemble $X(K)$ muni de la topologie \og forte\fg\ d\'eduite de la valeur absolue.
\begin{thm}\label{ThIntro}
Soit $(K,\vert\,.\,\vert)$ un corps valu\'e alg\'ebriquement clos, et soit $f:X\to Y$ un morphisme de $K$-sch\'emas de type fini. Les conditions suivantes sont \'equivalentes:
\begin{romlist}
\item\label{ThIntro1} $f$ est  \emph{universellement ouvert};
\item\label{ThIntro2} l'application continue  $f_{\Top}:X_{\Top}\to Y_{\Top}$ induite par $f$ est \emph{ouverte}.
\end{romlist}
\end{thm}
\begin{rems}\label{IntroRem}
\item\label{IntroRem2} Soit $f: X\to Y$ un morphisme localement de type fini de sch\'emas localement noeth\'eriens. Rappelons que $f$ est universellement ouvert dans les cas suivants:
\begin{romlist}
\item\label{IntroRem21} $f$ est plat (ou, plus g\'en\'eralement, il existe un $\cO_{X}$-module coh\'erent  et $f$-plat, de support $X$), cf. \cite{EGA4_II}, (2.4.6);
\item\label{IntroRem22} $Y$ est g\'eom\'etriquement unibranche et $f$ est \'equidimensionnel  (\og crit\`ere de Chevalley\fg\ \cite{EGA4_III}, (14.4.4)). On peut d'ailleurs, pour beaucoup de questions, se ramener \`a ce cas en utilisant par exemple (\cite{EGA4_III}, (14.4.9)).
\end{romlist}
\item\label{IntroRem3} Tout morphisme plat $f:X\to Y$ d'espaces \emph{analytiques complexes} est ouvert; ce r\'esultat semble d\^u \`a Douady (\cite{Dou68}, corollaire \`a la fin de l'article). On en trouvera une d\'e\-mons\-tra\-tion plus \'el\'ementaire dans  \cite{BaSta76} (V, Theorem 2.12), ou \cite{BaSta77vol2} (V, th\'eor\`eme 2.10), qui montre en outre (sans hypoth\`ese de platitude) que si $Y$ est localement irr\'eductible et $f$ \'equidimensionnel, alors $f$ est ouvert: on reconna\^{\i}t l'analogue de la condition \ref{IntroRem22} de \ref{IntroRem2}.
\item\label{IntroRem4} Tout morphisme plat d'espaces \emph{analytiques rigides} (quasi-compacts et quasi-s\'epar\'es) est ouvert (\cite{Bosch09}, cor. 7.2); la d\'emonstration utilise le cas d'un morphisme de sch\'emas via les mod\`eles formels.
\item\label{IntroRem5} Selon V. Berkovich (non publi\'e, communication priv\'ee), si $f:X\to Y$ est un morphisme \emph{sans bord} d'espaces analytiques non archim\'ediens, et si $\cF$ est un $\cO_{X}$-module coh\'erent et $f$-plat, alors la restriction de $f$ au support de $\cF$ est une application ouverte.
\item\label{IntroRem6} Vu les trois remarques pr\'ec\'edentes, la partie \og \ref{ThIntro1}$\,\Rightarrow\,$\ref{ThIntro2}\fg\  du  th\'eor\`eme \ref{ThIntro} est essentiellement connue %
dans le cas plat (et dans le cas complexe); %
les seuls m\'erites de la d\'emonstration pr\'esent\'ee ici sont %
(outre le cas ultram\'etrique non plat) %
de ne pas traiter s\'epar\'ement  les cas archim\'edien et non archim\'edien, et de ne pas recourir  aux espaces analytiques sous-jacents. Quant \`a la partie \og \ref{ThIntro2}$\,\Rightarrow\,$\ref{ThIntro1}\fg\ c'est un simple exercice, que je n'ai toutefois pas trouv\'e dans la litt\'erature.
\item Je remercie Vladimir Berkovich, Siegfried Bosch, Antoine Ducros, Francesco Polizzi et le rapporteur pour leurs remarques.
\end{rems}

\section{Rappels et pr\'eliminaires}\label{SecPrelim}
\subsection{Corps topologiques}\label{SsecCorpsTop}
Soit $K$ un corps topologique; on supposera toujours dans la suite que la topologie de $K$ n'est pas grossi\`ere, ce qui implique que $K$ est s\'epar\'e (d\'emonstration: $K$ admet un ouvert non vide et strict; par translation, il admet un ouvert non vide $U\subset K^*$; donc $K^*=\bigcup_{\lambda\in K^*}\lambda U$ est ouvert, donc $\{0\}$ est ferm\'e).

Pour tout $K$-sch\'ema de type fini $X$, on sait alors munir l'ensemble $X(K)$ d'une topologie dite \emph{forte} de telle sorte que, si l'on note $X_{\Top}$ l'espace topologique obtenu:
\begin{itemize}
\item si $X=\Aa^1_{K}$, la topologie forte sur $X_{\Top}=K$ est celle de $K$;
\item tout $K$-morphisme $f:X\to Y$ de $K$-sch\'emas de type fini induit une application \emph{continue} $f_{\Top}:X_{\Top}\to Y_{\Top}$, qui est de plus un plongement topologique ouvert (resp. ferm\'e) si $f$ est une immersion ouverte (resp. ferm\'ee);
\item si $X$ et $Y$ sont deux $K$-sch\'emas de type fini, la bijection canonique $(X\times_{K}Y)_{\Top}\to X_{\Top}\times Y_{\Top}$ est un hom\'eomorphisme.
\end{itemize}

Ces conditions entra\^{\i}nent encore que toute immersion induit un plongement topologique, et que la formation de $X_{\Top}$ commute aux produits fibr\'es.\medskip

Nous aurons \`a consid\'erer les propri\'et\'es suivantes d'un corps topologique $K$:\smallskip

\noindent (H) Pour tout $K$-morphisme \emph{\'etale} $f:X\to Y$ de $K$-sch\'emas de type fini, l'application induite  $f_{\Top}:X_{\Top}\to Y_{\Top}$ est \emph{ouverte}.\smallskip

\noindent (F) Pour tout $K$-morphisme \emph{fini} $f:X\to Y$ de $K$-sch\'emas de type fini, l'application induite  $f_{\Top}:X_{\Top}\to Y_{\Top}$ est \emph{ferm\'ee}.
\begin{subprop}\label{OuvDense} Soit $K$ un corps topologique \emph{non discret}. Soit $X$ un $K$-sch\'ema de type fini  \emph{irr\'eductible}, et soit $\Omega$ un ouvert non vide de $X_{\Top}$. Alors, dans chacun des cas suivants,  $\Omega$ est dense dans $X$ pour la topologie de Zariski:
\begin{romlist}
\item\label{OuvDense1} $X=\Aa^n_{K}$;
\item\label{OuvDense2} $K$ v\'erifie  \textup{(H)} et $X$ est lisse sur $K$;
\item\label{OuvDense3} $K$ v\'erifie  \textup{(H)} et est alg\'ebriquement clos.
\end{romlist}
\end{subprop}
\dem dans le cas \ref{OuvDense1}, l'hypoth\`ese que $K$ est non discret entra\^{\i}ne directement le r\'esultat si $n=1$ (tout ouvert non vide de $K$ est infini). Le cas $n\geq1$ s'en d\'eduit gr\^ace au fait bien connu suivant: si $I_{1},\dots,I_{n}$ sont des parties infinies de $K$, tout polyn\^ome $F\in K[T_{1},\dots,T_{n}]$ nul sur $\prod_{k=1}^n I_{k}$ est nul.\smallskip

Dans le cas \ref{OuvDense2}, soit $x$ un point de $\Omega$. Il existe un voisinage $U$ de $x$ dans $X$ et un morphisme \'etale $\pi:V\to\Aa^n_{K}$. D'apr\`es (H), $\pi(\Omega)$ est un ouvert non vide de $K^n$ donc est Zariski-dense dans $\Aa^n_{K}$, ce qui entra\^{\i}ne imm\'ediatement la conclusion.\smallskip

Dans le cas \ref{OuvDense3}, soit $x$ un point de $\Omega$. Il suffit de montrer que pour toute courbe irr\'eductible $C\subset X$ contenant $x$, l'ensemble $\Omega\cap C$ est infini. Nous sommes donc ramen\'es au cas o\`u $X$ est une courbe, et m\^eme  au cas d\'ej\`a trait\'e  d'une courbe lisse en passant \`a la normalis\'ee  (puisque $K$ est alg\'ebriquement clos, la normalisation induit une surjection sur les points rationnels).\qed

\subsection{Corps valu\'es}
Tout corps valu\'e  $(K,\vert\,.\,\vert)$ sera consid\'er\'e comme corps topologique, la topologie \'etant celle de la distance $(x,y)\mapsto\vert x-y\vert$ sur $K$. On sait qu'alors le compl\'et\'e $\wh{K}$ de $K$ est encore de fa\c{c}on naturelle un corps valu\'e, alg\'ebriquement clos si $K$ l'est.

\begin{subprop}\label{TopVal} Soit $(K,\vert\,.\,\vert)$ un corps valu\'e.
\begin{romlist}
\item\label{TopVal1} Si $K$ est s\'eparablement ferm\'e dans $\wh{K}$, il v\'erifie  la condition \textup{(H)} de \rref{SsecCorpsTop}. Plus pr\'e\-ci\-s\'e\-ment, tout $K$-morphisme \'etale $f:X\to Y$ de $K$-sch\'emas de type fini induit un \emph{hom\'eo\-mor\-phisme local} $f_{\Top}:X_{\Top}\to Y_{\Top}$.
\item\label{TopVal2} Si $K$ est alg\'ebriquement ferm\'e dans $\wh{K}$, il v\'erifie  la condition \textup{(F)} de \rref{SsecCorpsTop}.
\end{romlist} 
En particulier, tout corps valu\'e complet ou alg\'ebriquement clos v\'erifie \textup{(H)} et \textup{(F)}.
\end{subprop} 
Rappelons bri\`evement la d\'emonstration: on se ram\`ene, par des d\'evissages convenables, au cas o\`u $Y=\Spec A$ est l'espace affine de dimension $n$ sur $K$ et o\`u 
$$X=\begin{cases}
\Spec A[T,\frac{1}{F'(T)}]/(F(T)),\;F\in A[T] & \text{ (cas \ref{TopVal1})}\\
\Spec A[T]/(F(T)),\;F\in A[T]\text{ unitaire} & \text{ (cas \ref{TopVal2})}.
\end{cases}
$$
Lorsque $K$ est \emph{complet}, le cas \ref{TopVal1} est le th\'eor\`eme d'inversion locale classique, et le cas \ref{TopVal2} est la \og continuit\'e des racines\fg\ pour laquelle on renvoie par exemple \`a \cite{BGR84}, \S\ 3.4. 

Dans le cas g\'en\'eral, notons $\wh{f}_{\Top}:X(\wh{K})\to Y(\wh{K})$ l'application induite par $f$ sur les points \`a valeurs dans $\wh{K}$: on remarque alors que si $K$ est s\'eparablement (resp. al\-g\'e\-bri\-que\-ment) ferm\'e dans $\wh{K}$, on a  $\wh{f}_{\Top}^{-1}(Y(K))=X(K)$ dans le cas \ref{TopVal1} (resp. \ref{TopVal2}), de sorte que l'on est ramen\'e au cas complet.\qed

\begin{subrem}
La condition que $K$ soit s\'eparablement ferm\'e dans $\wh{K}$ \'equivaut au fait que la valeur absolue $\vert\,.\,\vert$ soit \emph{hens\'elienne}, c'est-\`a-dire admette un unique prolongement \`a toute extension alg\'ebrique de $K$. Dans le cas non archim\'edien, il revient au m\^eme de dire que l'anneau de la valuation associ\'ee (i.e. la boule unit\'e ferm\'ee de $\vert\,.\,\vert$) est hens\'elien. Voir par exemple  \cite{BAC}, chapitre 6, exercices du \S\ 8.
\end{subrem}

\begin{subrem}
Supposons la valeur absolue de $K$ associ\'ee \`a une valuation discr\`ete d'anneau $A$  hens\'elien. 

Si de plus $A$ est  \emph{excellent}. alors $K$ est m\^eme alg\'ebriquement ferm\'e dans $\wh{K}$, puisque ce dernier est une extension s\'eparable de $K$. Les conditions (H) et (F) sont donc v\'erifi\'ees dans ce cas. 

En revanche, il y a des exemples, en caract\'eristique $p>0$,  o\`u $\wh{K}$ contient un \'el\'ement  $\alpha\not\in K$ tel que $\alpha^p\in K$; dans ce cas, le morphisme de Frobenius (relatif) $f:\Aa^1_{K}\to\Aa^1_{K}$ est fini, mais $f_{\Top}$ n'est pas ferm\'ee puisque $\alpha^p$ est adh\'erent \`a $K^p$ sans appartenir \`a $K^p$.
\end{subrem}

\subsection{Applications et morphismes universellement ouverts}
Commen\c{c}ons par un \'enonc\'e topologique \'el\'ementaire. Nous dirons qu'une application continue $p:U\to V$ est un \emph{quotient topologique} si $V$ s'identifie \`a l'espace quotient de $U$ par la relation d'\'equivalence associ\'ee \`a $p$. C'est notamment le cas si $p$ est \emph{surjective et ouverte}, ou \emph{surjective et ferm\'ee}. %
\begin{subprop}\label{DescTop}
Consid\'erons un diagramme \emph{cart\'esien} d'applications continues
$$\xymatrix{%
U'\ar[d]^{\varphi'}\ar[r]^{p'} & U\ar[d]^{\varphi}\\
V'\ar[r]^{p} & V.\\
}$$
\begin{romlist}
\item\label{DescTop1} Si $\varphi$ est ouverte, alors $\varphi'$ est ouverte. \emph{(Autrement dit, les applications continues ouvertes sont universellement ouvertes).}
\item\label{DescTop2} Si $p$ est un quotient topologique et si $\varphi'$ est ouverte, alors $\varphi$ est ouverte.
\end{romlist}
\end{subprop}
\dem pour \ref{DescTop1}, on se ram\`ene \`a deux cas:
\begin{itemize}
\item $T$ est un espace topologique, $V'=V\times T$, et $p$ est la projection naturelle: alors  $U'=U\times T$ et l'assertion r\'esulte de la d\'efinition de la topologie produit.
\item $V'$ est un sous-espace de $V$ et $p$ est l'inclusion: alors $U'=\varphi^{-1}(V')$. Un ouvert $W'$ de $U'$ est de la forme $W\cap U'$, o\`u $W$ est un ouvert de $U$. L'assertion r\'esulte de la \og formule de projection\fg\ $\varphi(W\cap \varphi^{-1}(V'))=\varphi(W)\cap V'$.
\end{itemize}

Montrons \ref{DescTop2}: si $W$ est un ouvert de $U$, alors, le diagramme \'etant cart\'esien, on a $p^{-1}(\varphi(W))=\varphi'({p'}^{-1}(W))$. On  en d\'eduit que $p^{-1}(\varphi(W))$ est ouvert dans $V'$, d'o\`u la conclusion vu l'hypoth\`ese sur $p$.\qed
\smallskip

En ce qui concerne les morphismes universellement ouverts de sch\'emas, nous aurons besoin de la proposition suivante, extraite de  \cite{EGA4_III}, (14.5.9), qui affirme l'existence de \og quasi-sections\fg\ pour de tels morphismes:
\begin{subprop}\label{14.5.9} Soient $Y$ un sch\'ema localement noeth\'erien irr\'eductible, $f:X\to Y$ un morphisme localement de type fini, $y$ un point g\'eom\'etriquement unibranche (par exemple \emph{normal}) de $Y$. Les conditions suivantes sont \'equivalentes:
\begin{romlist}
\item \label{14.5.91} $f$ est universellement ouvert en tout point de $f^{-1}(y)$;
\item \label{14.5.92} pour tout point ferm\'e $x$ de $f^{-1}(y)$, il existe une partie localement ferm\'ee irr\'eductible $X'$ de $X$, contenant $x$, et telle que la restriction $X'\to Y$ de $f$ soit un morphisme quasi-fini et dominant.\qed
\end{romlist}
\end{subprop}

Enfin, rappelons une technique de localisation \'etale qui permet de \og rendre finis les morphismes quasi-finis\fg:

\begin{subthm}[\cite{EGA4_IV}, (18.12.1)]\label{LocEtale} Soient $f:X\to Y$ un morphisme localement de type fini, $y$ un point de $Y$, $x$ un point isol\'e de la fibre $f^{-1}(y)$. Alors il existe un morphisme \'etale $Y'\to Y$, un point $x'$ de $X'=X\times_{Y}Y'$ au-dessus de $x$, et un voisinage ouvert $V'$ de $x'$ dans $X'$ tel que le morphisme compos\'e $V'\inj X'\xrightarrow{\pr_{2}} Y'$ soit fini.\qed
\end{subthm}
\medskip

D'apr\`es \ref{TopVal}, le th\'eor\`eme \ref{ThIntro} r\'esulte des suivants, un peu plus pr\'ecis:
\begin{thm}\label{ThFacile} Soit $K$ un corps topologique \emph{alg\'ebriquement clos, non discret et  v\'erifiant \textup{(H)}}, et soit $f:X\to Y$ un morphisme de $K$-sch\'emas de type fini.  Si $f_{\Top}$  est ouverte, alors  $f$ est universellement ouvert.
\end{thm}
\begin{thm}\label{ThPpal} Soit $K$ un corps topologique \emph{alg\'ebriquement clos et v\'erifiant \textup{(H)} et  \textup{(F)}}, et soit $f:X\to Y$ un morphisme de $K$-sch\'emas de type fini.  Si $f$ est universellement ouvert, alors $f_{\Top}$  est ouverte.
\end{thm}
\noindent{\sl Preuve de \textup{\ref{ThFacile}}} (le th\'eor\`eme \ref{ThPpal} sera d\'emontr\'e au paragraphe  \ref{SecUnivOuv}).  D'apr\`es \ref{DescTop}\,\ref{DescTop1} il suffit de montrer, sous les  hypoth\`eses de \ref{ThFacile}, que $f$ est un morphisme ouvert et m\^eme, quitte \`a restreindre $X$, que $f(X)$ est ouvert dans $Y$. Comme $f(X)$ est constructible, il suffit de voir que pour tout ferm\'e irr\'eductible $Z$ de $Y$, l'ensemble $f(X)\cap Z$ est vide ou dense dans $Z$. On peut m\^eme supposer, par changement de base, que $Z=Y$ et notre assertion r\'esulte alors de \ref{OuvDense}\,\ref{OuvDense3}.\qed

\section{Preuve du th\'eor\`eme \ref{ThPpal}}\label{SecUnivOuv}
Dans tout ce qui suit, $K$ d\'esigne un corps topologique alg\'ebriquement clos. 
\subsection{Le cas fini et le cas quasi-fini}
Commen\c{c}ons par \'etablir \ref{ThPpal} lorsque le morphisme $f$ est fini, avec une hypoth\`ese un peu plus faible sur $K$ (toujours suppos\'e alg\'ebriquement clos):
\begin{subprop}\label{PropCasFini}
On suppose que $K$ v\'erifie la condition \textup{(F)}. Soit $f:X\to Y$ un morphisme \emph{fini et universellement ouvert} de $K$-sch\'emas de type fini. Alors l'application $f_{\Top}$ est ouverte.
\end{subprop}
\dem la conclusion est imm\'ediate lorsque $X$ est r\'eunion des images d'une famille de sections de $f$; nous allons nous ramener \`a ce cas par changement de base.

Fixons un entier $d$ majorant le cardinal de toutes les fibres de $f_{\Top}$. Con\-si\-d\'e\-rons le diagramme suivant de $Y$-sch\'emas finis, o\`u les produits et puissances sont \emph{fibr\'es sur $Y$}:
$$\xymatrix{%
X\times X^d \ar[r]^<(0.4){\mathrm{pr}_{1}} \ar[d]^{f_{1}=\mathrm{pr}_{2}}
&X \ar[d]^{f}\\
X^d\ar[r] \ar@/^/@<2pt>[u]\ar@/^/@<5pt>[u]\ar@/^/@<8pt>^{s_{1},\dots,s_{d}}[u]
& Y
}$$
et o\`u les sections $s_{1},\dots, s_{d}$ de la projection $f_{1}$ sont d\'efinies par 
$$s_{i}(x_{1},\dots,x_{d})=(x_{i},x_{1},\dots,x_{d})\quad (1\leq i\leq d).$$
La r\'eunion $S$ des images des $s_{i}$ est un ferm\'e de $X\times X^d$; comme $f_{1}$ est ouvert vu l'hypoth\`ese sur $f$, il en r\'e\-sulte que
$$Y':=\{t\in X^d \mid f_{1}^{-1}(t)\subset S\}$$
est un \emph{ferm\'e} de $X^d$ (son compl\'ementaire est l'image par $f_{1}$ du compl\'ementaire de $S$). Vu le choix de $d$, $Y'$ s'envoie surjectivement sur $Y$: si $y\in Y(K)$, alors $f^{-1}(y)$ est de la forme $\{x_{1},\dots,x_{d}\}$ (avec des r\'ep\'etitions \'eventuelles) et le point $(x_{1},\dots,x_{d})$ au-dessus de $y$ appartient \`a $Y'$. 
Par restriction \`a $Y'$, on obtient donc un diagramme cart\'esien
$$\xymatrix{%
X'\,\ar@{^{(}->}@<-2pt>[r] \ar[d]^{f'}&X\times X^d \ar@<-2pt>[r]^<(0.4){\mathrm{pr}_{1}} \ar[d]^{f_{1}}
&X \ar[d]^{f}\\
Y'\,\ar@{^{(}->}@<-2pt>[r] \ar @/_/@<-5pt>[rr]_{p\text{ (surjectif et ferm\'e)}}&X^d\ar@<-2pt>[r] 
& Y
}$$
o\`u, par construction, $X'$ est r\'eunion des images des sections induites par les $s_{i}$, de sorte que $f'_{\Top}$ est ouverte. L'application  $p_{\Top}$ est surjective ($K$ est alg\'ebriquement clos) et ferm\'ee ($K$ v\'erifie (F)) donc est un quotient topologique. Donc $f_{\Top}$ est ouverte en vertu de \ref{DescTop}\,\ref{DescTop2} et du fait que le foncteur $Z\mapsto Z_{\Top}$ commute aux produits fibr\'es finis.\qed \smallskip\\
Par localisation \'etale, on en d\'eduit le cas \emph{quasi-fini} de \ref{ThPpal}:
\begin{subcor}\label{CorCasQFini}
Sous les hypoth\`eses du th\'eor\`eme \rref{ThPpal}, on suppose en outre que $f$ est \emph{quasi-fini}. Alors $f_{\Top}$ est ouverte.
\end{subcor}
\dem fixons un point $x$ de $X_{\Top}$, d'image $y$ dans $Y_{\Top}$, et montrons que $f_{\Top}$ est ouverte au point $x$. Si $p:Y'\to Y$ est un morphisme \'etale, $y'$ un point de $Y'_{\Top}$ au-dessus de $y$, et $x'$ le point $(x,y')$ de $X\times_{Y}Y'$, il suffit de montrer l'\'enonc\'e analogue pour $x'$ et $\mathrm{pr}_{2}: X\times_{Y}Y'\to Y'$,  compte tenu de \ref{DescTop}\,\ref{DescTop2} et de l'hypoth\`ese (H) faite sur $K$.

Or, d'apr\`es \ref{LocEtale}, il existe un tel couple $(Y',y')$ et un ouvert $Z'\subset X\times_{Y}Y'$ contenant le point $x'$ et \emph{fini} sur $Y'$: il suffit alors d'appliquer \ref{PropCasFini}.\qed
\subsection{R\'eduction au cas d'un morphisme quasi-fini, de but normal}
Commen\c{c}ons par ramener la d\'emonstration de \ref{ThPpal} au cas o\`u $Y$ est \emph{normal}. Soit $p: Y'\to Y$ le normalis\'e de $Y_{\mathrm{red}}$ et soit $f':X'\to Y'$ d\'eduit de $f$ par le changement de base $p$: alors $p_{\Top}$ est ferm\'ee et surjective donc est un quotient topologique. Ainsi il suffit, d'apr\`es \ref{DescTop}, de d\'emontrer que $f'_{\Top}$ est ouverte. 

Il nous reste donc \`a traiter le cas d'un morphisme $f:X\to Y$ o\`u $Y$ est normal.  

\subsection{Fin de la d\'emonstration du th\'eor\`eme} 
Sous les hypoth\`eses de \ref{ThPpal}, et en supposant $Y$ normal,  soit $x\in X(K)$ un point quelconque: appliquant \ref{14.5.9}, on trouve un sous-sch\'ema localement ferm\'e irr\'eductible $X'$ de $X$ contenant $x$, quasi-fini et dominant sur $Y$. En outre, l'implication \ref{14.5.92}\,$\Rightarrow$\,\ref{14.5.91} de \ref{14.5.9} montre que le morphisme induit $X'\to Y$ est encore universellement ouvert.  D'apr\`es \ref{CorCasQFini}, l'application induite $X'_{\Top}\to Y_{\Top}$ est ouverte, donc $f_{\Top}$ est ouverte au point $x$.\qed

\bibliography{biblio}
\end{document}